\documentclass[12pt,a4paper]{article}
\RequirePackage{etex}
\usepackage{amsfonts, amsmath, amssymb,amsthm,graphics,mathrsfs,verbatim, CJK}
\usepackage{makeidx}
\usepackage{palatino}
\usepackage{titletoc}
\usepackage{tocloft}
\usepackage{titlesec}
\usepackage{stmaryrd,color}
\usepackage[colorlinks=true]{hyperref}
\usepackage{graphicx}
\usepackage{pgfplots}
\usepackage{tikz}
\usepackage{pstricks,pst-node,multido,ifthen,calc}

\titleformat{\section}{\Large}{}{0.5em}{}
\makeindex

\normalsize
\pgfplotsset{compat=1.14}

\def\b{\beta}\def\D{\Delta}
\def\eset{\emptyset}

\def\ol{\overline}

\def\nin{\noindent}

\def\seq{\subseteq}

\def\fs{\footnotesize}

\def\vs{\vspace*}

\newtheorem{mthm}{Theorem}[section]
\newtheorem{mylem}[mthm]{Lemma}
\newtheorem{myprn}[mthm]{Proposition}
\newtheorem{mycor}[mthm]{Corollary}
\newtheorem{mydef}[mthm]{Definition}
\newtheorem{myrem}[mthm]{Remark}
\newtheorem{mycon}[mthm]{Construction}
\newtheorem{myeg} [mthm]{Example}
\newtheorem{myque} [mthm]{Question}
\newtheorem{myalg} [mthm]{Algorithm}
\newtheorem{myconj} [mthm]{Conjecture}

\newenvironment{prof}{\noindent $Proof.$ \rm}{\hfill $\Box$}

\def\b{\beta}\def\D{\Delta}\def\fs{\footnotesize}

\def\seq{\subseteq}

\def\ol{\overline}

\def\vs{\vspace*}

\def\varnothing{\emptyset}

\def \nin {\noindent}

\def \Lemma #1 {\vs{2mm}\nin {\bf Lemma #1} }
\def \Prop #1 {\vs{2mm}\nin {\bf Proposition #1} }
\def \Th #1 {\vs{2mm}\nin {\bf Theorem #1} }
\def \Cor #1 {\vs{2mm}\nin {\bf Corollary #1} }

\def \Def #1 {\vs{2mm}\nin {\bf Definition #1} }
\def \part #1 {\hfil\break\hglue 12pt {\rm (#1)~}}

\def\fs{\footnotesize}

\def\vd{vertex decomposable}
\def\cwg{clique-whiskered graph}

\def\mcwg{multiple clique-whiskered graph}
\def\mclwg{multiple cluster-whiskered graph}
\def\vcp{vertex clique-partition}

\def\Epi{{cc}}
\def\Epi1m{cc}
\def\epi2m{cc}

\begin{document}
\title{
\bf\LARGE  A construction of sequentially Cohen-Macaulay graphs \thanks{This research was supported by the National Natural
Science Foundation of China (Grant No. 11271250). } }
\author{{A-Ming Liu\thanks{  aming8809@163.com} and Tongsuo Wu\thanks{Corresponding author. tswu@sjtu.edu.cn}}\\
 {\small School of Mathematical Sciences,  Shanghai Jiao Tong University
}
}

\date{}
\baselineskip=16pt \maketitle

\begin{center}
\begin{minipage}{12cm}

\vs{3mm}\nin{\small\bf Abstract.} {\fs For every simple graph $G$, a class of multiple clique cluster-whiskered graphs $G^{md}$ is introduced, and it is shown that all graphs $G^{md}$ are \vd, thus the independence simplicial complex ${\rm Ind}\,G^{md}$ is  sequentially Cohen-Macaulay; the properties of the graphs $G^{md}$ and the clique-whiskered graph $G^\pi$ are studied, including the enumeration of facets of the complex ${\rm Ind}\, G^{\pi}$ and the calculation of Betti numbers of the cover ideal $I_c(G^{md})$ for some graph $G^{md}$. }

\vs{3mm}\nin {\small Key Words and Phrases:} {\small  Multiple clique cluster-whiskered graph; vertex decomposable graph; sequentially Cohen-Macaulay; facet ideal}

\vs{3mm}\nin {\small 2010 AMS Classification:} {\small   Primary: 13H10; 05E45; Secondary: 13F55; 05E40.}

\end{minipage}
\end{center}

\section{1. Introduction}

Vertex decomposable simplicial complexes and graphs are important in combinatorial topology  and combinatorial commutative algebra, since they provide examples of shellable simplicial complexes, while  vertex decomposable simplicial complexes and graphs are further sequentially Cohen-Macaulay, see \cite{Herzog and Hibi 1999} for a characterization of sequentially Cohen-Macaulay complexes.

For each graph $G$, Villarreal \cite{Villareal1990} constructed a full whiskered graph $G^W,$  i.e., for each vertex $v$ of $G$, an end vertex is attached to $v$ (the induced subgraph $K_2$ is called a {\it whisker}) and,  it is proved that the new graph $G^W$ is Cohen-Macaulay. Later, Dochermann and Engste\"om \cite{ADAE} proved that $G^W$ is vertex decomposable.
Cook and Nagel \cite{Cook2011} generalized the whiskered idea and constructed the vertex \cwg \,$G^\pi$ and proved, among other things, that $G^\pi$ is unmixed and \vd; that the $f$-vector of $G$ is the $h$-vector of $G^\pi$.
Biermann and Van Tuyl \cite{BT2013} defined the concept of $s$-coloring $\chi$ on a simplicial complex $\D$ to construct a balanced pure \vd \,simplicial complex $\D_{\chi}$,  generalizing the construction of ${\rm Ind}\,G^\pi$ and establishing some similar results of \cite{Cook2011}  to the general setting of simplicial complexes. See also \cite{2014BC} for edge-clique-whisker graphs and related \vd\, property.

 Recall that Francisco and  Tuyl in 2007 \cite{HA} proved that chordal graphs are sequentially Cohen-Macaulay; later it was further proved   that chordal graphs are shellable in 2008 by Tuyl and Villarreal \cite{VanVillareal}, and that chordal graphs are \vd\, in 2009 by Woodroofe \cite{Woodroofe}. Francisco and H\'a \cite{FH2008} studied the sequentially Cohen-Macaulay property of  graphs obtained after partial whiskers added. Biermann, Francisco, H\'a  and Tuyl in \cite{BCFHT2015} continued the works of  \cite{BT2013} and constructed a partial coloring on simplicial complexes, thus providing a general way for constructing classes of sequentially Cohen-Macaulay simplicial complexes.

 Recall also that Francisco, H\'a and Tuyl in \cite{FHT2009} have had a systematic study on the Betti splitting of monomial ideals; a particular case of Betti splitting, called {\it variable splitting }, was studied in \cite{Moradi2} for monomial ideals, and it is proved that a simplicial complex $\D$ is \vd\, if and only if the ideal $I_{\D^\vee}$ admits variable splitting property; such ideals $I_{\D^\vee}$ further has Betti splitting property, thus the graded Betti numbers of $I(\ol{\D})$ can be calculated recursively if $\D$ is \vd.

This work is a continued study of clique-whiskered graphs. In section 2, we recall some background used in the paper; in section 3, we study an exchange property on the facets of the complex $\D={\rm Ind}\, G^\pi$ and enumerate the facets of $\D$; in section 4, we construct a class of multiple clique cluster-whiskered graphs $G^{md}$, which includes $G^\pi$ as particular cases and, show that it is nonpure \vd, and we study some invariants of $G^{md}$.

\section{2. Preliminaries}

In this section, we introduce some relevant definitions and results on simplicial complexes, clutters, and minimal free resolutions. We use \cite{Herzog and Hibi, Villareal} as basic references.

Recall that a \emph{simplicial complex} $\D$ on the vertex set $[n]=\{1,2,\ldots, n \}$ is a collection of subsets of $[n]$ such that if $F\in \D$ and $E \subseteq F$, then $E\in \D$. Each $F$ in $\D$ is called a \emph{face}, and a \emph{facet} $F$ is a maximal face with respect to inclusion. A simplicial complex $\D$ is called \emph{pure} if all facets have the same cardinality. Denote by $\mathcal{F}(\D)$  the set of all  facets of $\D$ and, if $\mathcal{F}(\D)=\{F_1, F_2, \ldots, F_t\}$, then $\D$ can be denoted by $\langle F_1, F_2, \ldots, F_t \rangle$. A \emph{nonface} of $\D$ is a subset $F$ of $[n]$ with $F\notin\D$. Let $\mathcal{N}(F)$ be the set of minimal nonfaces of $\D$ with respect to inclusion. The \emph{dimension} of a face $F$ is $dim \,F=|F|-1$ and, the dimension of $\D$ is defined by $ dim\,\D=max\{dim F \mid F\in \D \}$. Recall that the Alexander dual $\D^\vee$ of a simplicial complex $\D$ is defined by
$\D^\vee=\{[n]\smallsetminus F \mid F\notin \D\}.$

For a face $H$ of a simplicial complex $\D$, recall the following notation
$$\D\smallsetminus H=:\{F\in \D \mid H\cap F=\eset\},$$
$${\rm lk}_{\D}(H)=:\{F\in \D \mid H\cap F=\eset,\, F\cup H \in \D\}.$$
Note that whenever $H=\{x\}$ is a vertex, the notations are usually written as $\D\smallsetminus x$ and, ${\rm lk}_{\D}(x)$ respectively.
Recall  the following  concept of a {\it vertex decomposable} simplicial complex, which is introduced by Provan and Billera \cite{PB} in pure case and, is extended to the nonpure case by Bj\"{o}rner and Wachs \cite{BW}:

\begin{mydef}  A simplicial complex $\D$  with vertex set $[\,n\,]$ is called {\it vertex decomposable}, if either $\D$ is a simplex, or there is a vertex $x$ of $\D$ such that the following requirements are satisfied:

$(\alpha)$ Both $\D\smallsetminus x$ and $lk_{\D}(x)$ are  vertex decomposable;

$(\beta)$ No facet of $lk_{\D}(x)$ is a facet of $\D\smallsetminus x$, or equivalently,
$$\D\smallsetminus x=\langle\,\{F \mid x\not\in F\in \mathcal{F}(\D)\}\,\rangle.$$

\end{mydef}

A vertex $x$ satisfying conditions $(\alpha)$ and $(\beta)$ is called a {\it shedding vertex} of $\D$.  If $x$ only satisfies condition $(\beta)$, then we call it a {\it weak shedding vertex}.

A translation in the language of graphs is restated in the following

\begin{mydef}\label{VD graphversion}
A graph $G$ is called vertex decomposable if either it has no edges, or else has some vertex\index{graph!vertex decomposable}
$x$ such that we have as follows:

$(1)$ Both $G \smallsetminus N_G[x]$	 and $G\smallsetminus x$ are vertex decomposable, where $N_G[x]$ is the union of the neighbourhood $N_G(x)$ together with $\{x\}$;

$(2)$ For every independent set $S$ in $G \smallsetminus N_G[x]$, there exists some $y \in N_G(x)$ such that
$S\cup \{y\}$ is independent in $G\smallsetminus x$.
\end{mydef}

A vertex $x$ with the properties is called a {\it shedding vertex} of $G$.

It is well known that a {\it vertex cover} of a graph $G$ is a subset $C$ of the vertex set $V(G)$ such that $C\cap \{i,\,j\}\not=\varnothing$ holds for all $\{i,\,j\}\in E(G).$ A vertex cover $S$ is called minimal if no proper subset of $S$ is a vertex cover. Clearly, a subset $C$ of $V(G)$ is a minimal vertex cover of $G$ if and only if $V(G)\smallsetminus C$ is a maximal independent vertex subset of $G$.  Recall that a graph $G$ is said to be {\it unmixed} (or alternatively,{\it well-covered}), if all minimal vertex covers of $G$ have the same cardinality. It is known that a Cohen-Macaulay graph is always unmixed, while a graph $G$ is unmixed if and only if  the simplicial complex ${\rm Ind}\,G$ is pure.
The following result is well known:
\begin{myprn} \label{unmixed graph}A graph $G$ is called {\bf unmixed} or {\bf nonpure} if it satisfies one of the following equivalent conditions:\index{graph!unmixed characterization}

$(a)$ All the minimal vertex covers of the graph $G$ have a same cardinality.

$(b)$ The clique complex $\D_G$ of $\ol{G}$ (i.e., the independence simplicial complex ${\rm Ind}\,G$ of $G$) is pure, i.e., all the maximal independent vertex subsets have a same cardinality.

$(c)$ The edge ideal $I(G)$ of $G$ is height unmixed.
\end{myprn}

%
%
%
%

Recall the following

\begin{mydef} A simplicial complex $\D$ is called {\it shellable}, if its facets can be linearly ordered by
$F_1, F_2, \ldots, F_t,$ such that for any $2\le j\le t$ and $1\le i<j$, there exists an integer $1\le k<j$ satisfying

$(1)$ $|F_j\smallsetminus F_k|=1$ and,

$(2)$ $F_j\smallsetminus F_k \subseteq F_j\smallsetminus F_i$.

An order of the facets satisfying the conditions is called a {\it shelling order}.
\end{mydef}



For pure simplicial complexes, recall the following established implications:
$$vertex\,\,decomposable \,\, and \,\, pure\Longrightarrow pure \,shellable\Longrightarrow Cohen-Macaulay.$$
For nonpure simplicial complexes, it is known that \vd\, implies nonpure shellable, while nonpure shellable implies sequentially Cohen-Macaulay.

For  a face $F$ of a simplicial complex $\D$, let $X_F=\prod \limits_{i\in F}x_i$. Recall that the Stanley-Reisner ideal of $\D$ is defined by $$I_\D=\langle X_F \mid F\in \mathcal{N}(\D)\rangle.$$ Furthermore, if $\mathcal{F}(\D)=\{ F_1,F_2,\ldots,F_m\}$, then the facet ideal $I(\D)$ of $\D$ is defined by $I(\D)=\langle X_{F_1},X_{F_2},,\ldots,X_{F_m}\rangle.$ Clearly, $I_{\D^{\vee}}=I(\overline{\D})$ holds, where $\mathcal{F}(\overline{\D})=\{[n]\smallsetminus F\mid F\in \mathcal{F}(\D) \}.$

For a graph $G$, note that $G$ is said to be vertex decomposable (shellable, strongly shellable, (sequentially) Cohen-Macaulay, respectively), if the simplicial complex ${\rm Ind}\,G$ is vertex decomposable (shellable, strongly shellable, (sequentially) Cohen-Macaulay, respectively). Assume further $G$ has vertex set $\{x_1,\ldots,x_n\}$ with edge set $E(G)$, and let $S=k[x_1,x_2,\ldots,x_n]$ be the polynomial ring with $n$ variables $x_i$ over a field $k$. Recall that the {\it edge ideal} associated to $G$ is defined to be the ideal
$I(G)=\langle\, x_ix_j \mid \{x_i,x_j\}\in E(G)\,\rangle$
of $R$, and it is clear that $I(G)$ is identical with the Stanley-Reisner ideal $I_\D$ of the simplicial complex $\D={\rm Ind}\,G$. The edge ideal of a graph was introduced by Villarreal in \cite{Villareal}, and the Cohen-Macaulay property of the edge ideals was studied.

H\`{a} and Van Tuyl in \cite{HA} studied the resolutions of edge ideals, finding some recursive methods for computing Betti numbers of triangular hypergraphs; Fr\"oberg \cite{Froberg} gave a characterization on a property of the edge ideal $I(G)$, namely, {\it $I(G)$ has a linear resolution if and only if the complement graph $\ol{G}$ is chordal}.

%
%
%
%

\section{3. The maximal facets of the complex ${\rm Ind}\,G^\pi$ }

Recall that a clique of a graph is an induced complete subgraph. Recall from \cite{Cook2011} that a {\it vertex clique-partition} $\pi$ of a graph $G$ is a partition $\pi=\{W_1,\ldots,W_d\}$
of $V(G)$ such that each subgraph induced on $W_i$ is a nonempty clique. Note that $W_i$ in \cite{Cook2011} are allowed to be empty, and this is convenient for considering deletion and link when isolated vertices may occur. For the sake of computing Betti numbers, it is assumed that each $W_i$ in this paper is nonempty.

Let $G$ be any graph, let $\pi=\{W_1,\ldots,W_d\}$ be a vertex clique partition of $G$, and let $W=\{s_1',\ldots, s_d'\}$ be the vertices newly added in constructing the vertex clique-whiskered graph $G^\pi$, where each $s_i'$ is adjacent to all vertices of $W_i$. Set $\mathcal P(G^\pi)=\mathcal F({\rm Ind}\, G^\pi)$
 and, we consider enumeration of the facet set $\mathcal P(G^\pi)$; note that  the facets of ${\rm Ind}\,G^\pi$ have a kind of nice {\it exchange} property away from the  facet $W;$ we can define a pre-partial order $<$ in $\mathcal P(G^\pi)$, in the following way: in
$$F_1\leq F_2\,\, {\rm if\,\, and \,\,only\,\,if}\,\,F_1\smallsetminus W\seq F_2\smallsetminus W.$$
Note that this really defines a partial order in $\mathcal P(G^\pi)$, since $F_1\smallsetminus W=F_2\smallsetminus W$ does imply $F_1=F_2$, where $W$ is the least element of the poset $(\mathcal P(G^\pi),\leq )$. We have

\begin{myprn} \label{PowerSet} For any  maximal element $F$
of the poset $(\mathcal P(G^\pi),\leq)$, let $r=|F\smallsetminus W|$, where $W$ is the set of added vertexes in $G^\pi$ from $G$. Then we have

$(1)$ There exist  totally $r!$ maximal chains in $\mathcal P(G^\pi)$ from $W$ to $F$.

$(2)$ The number of facets of $G^\pi$ between $W$ and $F$ in $\mathcal P(G^\pi)$ is
 $2^r.$
\end{myprn}
\begin{prof} Let $F\smallsetminus W=\{v_1,\ldots,v_r\}$. Then each $G\in[W,F]$ is determined by the subset $G\smallsetminus W$ of $\{v_1,\ldots,v_r\}$. Assume that $G_1\smallsetminus W=\{v_1\}$, then the covers of $G_1$ in the poset $(\mathcal P(G^\pi),\leq)$ are determined by the subsets $\{v_1,v_i\}$ ($2\le i\le r$); assume that $G_2\smallsetminus W=\{v_2\}$, then the irredundant covers of $G_2$ in the poset $(\mathcal P(G^\pi),\leq)$ are determined by the subsets $\{v_2,v_i\}$ ($3\le i\le r$); and this proves the second conclusion.

The first conclusion follows in a similar way.
\end{prof}

\vs{3mm}Let $F_1,\ldots,F_s$ be the maximal elements of the poset, then clearly
$$\mathcal F({\rm Ind}\,G^\pi)=(\mathcal P(G^\pi),\leq)=\cup_{i=1}^s[W, F_i]$$
holds. Note that
 $F_i\smallsetminus W$ are all the maximal independent vertex set of the graph $G$, thus to count the number of maximal facets of the simplicial complex ${\rm Ind}\,G^\pi$, it is only necessary to count the number of the faces of the simplicial complex $\cup_{i=1}^s [\emptyset, (F_i\smallsetminus W)]$, and this depends on the mutual intersection property of the sets $F_1\smallsetminus W,\ldots,F_s\smallsetminus W$.

Note that in \cite[Proposition 3.8]{Cook2011}, it is proved that the $h$-vector of ${\rm Ind} \,G^\pi$ is identical with the $f$-vector of $G$.

We end this section with the following example:

 \begin{myeg}
 {\rm
 For the cycle graph  $G=C_{2n}$, let $G^\pi =: G^E$ be the clique-whiskered graph after adding $n$ ears. Clearly, $C_{6}$ has the following five maximal independent vertex sets:
 $$F_1=\{1,3,5\},F_2=\{2,4,6,\},F_3=\{1,4 \},F_4=\{ 2,5\},F_5=\{3,6 \};$$
thus the number of faces of $\cup_{i=1}^5[\emptyset,F_i]$ is
 $$2\times 2^3+3\times (2^2-2)-4\times 1=18.$$
By Proposition \ref{PowerSet}, the vertex independence complex of the \cwg\, $C_6^E$  has 18 facets.
}
 \end{myeg}


\section{4. Multiple clique cluster-partitioned graphs $G^{md}$ and $G^{mc}$ }

 Let $\pi=\{W_1,\ldots,W_d\}$ be a vertex clique-partition of a graph $G$, where each $W_i$ is not empty; let further $V(G)$ be a disjoint union of nonempty subsets $U_1,\ldots,U_s$ with the following additional
properties:
\begin{center}
\begin{minipage}{14cm}
$(1)$ For each $i$ ($1\le i\le d$) and $j$ ($1\le j\le s$), either $W_i\seq U_j$ or $W_i\cap U_j=\eset$.

$(2)$ For distinct $i,j$, if $W_i\cup W_j\seq U_k$, then the edge subset $E(W_i, W_j)$ is empty in $G$, where $E(W_i, W_j)=\{\{v_i,v_j\}\mid \{v_i,v_j\}\in E(G), v_i\in W_i, v_j\in W_j\}.$
\end{minipage}
\end{center}

\noindent Such a partition $e_\pi=\{U_1,\ldots, U_s\}$ is called a {\it clique cluster-partition} of $V(G)$ based on the vertex clique-partition $\pi$ of $G$. Note that in a clique cluster-partition $e_\pi=\{U_1,\ldots U_s\}$ based on a vertex clique-partition $\pi=\{W_1,\ldots,W_d\}$ of $V(G)$, $s\le t$ clearly holds true, while each  $U_i$  is either identical with some $W_j$ or,
is a union of more than one cliques from $\pi$ which are mutually independent in the graph $G$.

Let $e_\pi=\{ U_1,\ldots U_r \}$ be a clique cluster-partition of $V(G)$ based on the vertex clique-partition $\pi=\{W_1,\ldots,W_d\}$ of a graph $G$. Assume further that $U_j$ contains more than one $W_l$ for $j$ with $1\le j\le r$ and, for each $i$ with $r+1\le i\le s$, $U_i$ is identical with some $W_{k_i}$. Now we expand the graph $G$ to a new and larger graph $G^{md}$ in the following

\vs{2mm}{\bf Construction:} Let $A_1,\ldots,A_d$ be mutually disjoint nonempty finite sets with $A_i\cap V(G)=\eset\, (\forall i);$
let $B_{1},\ldots,B_r$ be mutually disjoint nonempty finite sets such that $$B_i\cap A_j=\eset,\,B_i\cap V(G)=\eset\,(\forall i,j).$$ Set
$$V(G^{md})=V(G)\cup A_1\cup\cdots\cup A_d\cup B_{1}\cup\cdots \cup B_r,$$
and the edge set $E(G^{md})$ are established in accordance with the following three rules:

$(1)$ $E(G)\seq E(G^{md})$ holds.

$(2)$ For each $a_i\in A_i,$ $W_i\cup \{a_i\}$ induces a clique
and, each vertex of $B_j$ is adjacent to all vertices of $U_j$, in $G^{md}$.

$(3)$ For each $j$, the subgraph of $G^{md}$ induced on $B_j$ is \vd; for each $i$, the subgraph of $G^{md}$ induced on $A_i$ is \vd.

\vs{2mm}The graph $G^{md}$ is called a {\it multiple clique cluster-whiskered graph } {\bf of type $(d,r)$}, abbreviated as a {\it mcc-whiskered} graph. A particular case is when the subgraphs of $G^{md}$ induced on each $B_j$ and on each $A_i$ are empty, i.e., $E_{G^{md}}(B_j)=\eset=E(A_i)$; in the case, the graph will be denoted as $G^{mc}$; if further $|A_i|=1=|B_j|$ hold for all possible $i,j$, the constructed graph $G^{mc}$ is denoted as $G^{cc}$ and is called a  {\it clique cluster-whiskered graph}, abbreviated as a {\it cc-whiskered graph}.
Note that if $d=s$ and, $U_i=W_i$ holds for all $i$, then
$G^{mc}$ is the multiple generalization of the clique whiskered graph $G^\pi$; if further $|A_i|=1$ holds true for each $i$, then $G^{cc}=G^{\pi}$, which is constructed and studied in \cite{Cook2011}.

We illustrate our construction $G^{mc}$ with the following example, where $G$ is taken to be a very simple graph $L_6$:

\begin{myeg}\label{mccExample}
{\rm
Let $G$ be the line graph of length 5:

\begin{center}
\begin{tikzpicture}[scale=0.9]

\draw(1,1)--(2,1)--(3,1)--(4,1)--(5,1)--(6,1);
\draw(1,1.25)node {$v_1$};
\draw(2,1.25)node {$v_2$};
\draw(3,1.25)node {$v_3$};
\draw(4,1.25)node {$v_4$};
\draw(5,1.25)node {$v_5$};
\draw(6,1.25)node {$v_6$};

\draw[fill](1,1)circle [radius=0.05];
\draw[fill](2,1)circle [radius=0.05];
\draw[fill](3,1)circle [radius=0.05];
\draw[fill](4,1)circle [radius=0.05];
\draw[fill](5,1)circle [radius=0.05];
\draw[fill](6,1)circle [radius=0.05];

\draw(3.5,0)node{\scriptsize Figure 4.1. The graph $G=L_6$};
\end{tikzpicture}
\end{center}
$(1)$ Consider the trivial clique decomposition  $\pi_1: V(G)=\cup\{v_i\}$ and, a clique cluster decomposition $e_{\pi_1}:V(G)=\{v_1,v_3,v_5\}\cup\{v_2,v_4,v_6\}$ based on $\pi_1$. Then we have the corresponding cc-whiskered graph $G^{cc}$:

\begin{center}
\begin{tikzpicture}[scale=0.9]

\draw(1,1)--(2,1)--(3,1)--(4,1)--(5,1)--(6,1);
\draw(1,0)--(1,1)--(3,2);
\draw(2,0)--(2,1)--(5,2);
\draw(3,0)--(3,1)--(3,2);
\draw(4,0)--(4,1)--(5,2);
\draw(5,0)--(5,1)--(3,2);;
\draw(6,0)--(6,1)--(5,2);

\draw(1,1.25)node {$v_1$};
\draw(2,1.25)node {$v_2$};
\draw(3.35,1.25)node {$v_3$};
\draw(3.9,1.25)node {$v_4$};
\draw(5,1.25)node {$v_5$};
\draw(6,1.25)node {$v_6$};

\draw(3,2.25)node {$b_1$};
\draw(5,2.25)node {$b_2$};

\draw(1,-0.25)node {$a_1$};
\draw(2,-0.25)node {$a_2$};
\draw(3,-0.25)node {$a_3$};
\draw(4,-0.25)node {$a_4$};
\draw(5,-0.25)node {$a_5$};
\draw(6,-0.25)node {$a_6$};

\draw[fill](1,1)circle [radius=0.05];
\draw[fill](2,1)circle [radius=0.05];
\draw[fill](3,1)circle [radius=0.05];
\draw[fill](4,1)circle [radius=0.05];
\draw[fill](5,1)circle [radius=0.05];
\draw[fill](6,1)circle [radius=0.05];

\draw[fill](1,0)circle [radius=0.05];
\draw[fill](2,0)circle [radius=0.05];
\draw[fill](3,0)circle [radius=0.05];
\draw[fill](4,0)circle [radius=0.05];
\draw[fill](5,0)circle [radius=0.05];
\draw[fill](6,0)circle [radius=0.05];

\draw[fill](3,2)circle [radius=0.05];
\draw[fill](5,2)circle [radius=0.05];

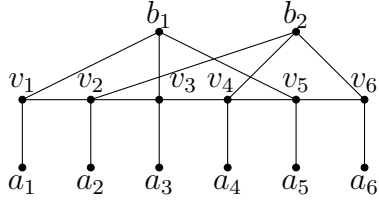
\draw(3.7,-1.5)node{\scriptsize Figure 4.2. The cc-whiskered graph $G^{cc}$ based on $e_{\pi_1}$};
\end{tikzpicture}
\end{center}

$(2)$ Consider the clique decomposition  $\pi_2: V(G)=\{v_3,v_4\}\cup\{v_1,v_2\}\cup \{v_5,v_6\}$ and, a clique cluster decomposition
$e_{\pi_2}:V(G)=\{v_1,v_2; v_5,v_6\}\cup\{v_3,v_4\}$ based on $\pi_2$. Let
$$A_1=\{a_{11},a_{12}\}, A_2=\{a_{21},a_{22},a_{2,3}\},A_3=\{a_3\}, B_1=\{b_{11},b_{12}\}.$$
Then we have the corresponding {\it mcc-whiskered graph} $G^{mc}$:

\begin{center}
\begin{tikzpicture}[scale=0.9]

\draw(1,1)--(2,1)--(3,1)--(4,1)--(5,1)--(6,1);

\draw(0.8,1.25)node {$v_1$};
\draw(2.2,1.25)node {$v_2$};
\draw(2.8,1.25)node {$v_3$};
\draw(4.2,1.25)node {$v_4$};
\draw(4.8,1.25)node {$v_5$};
\draw(6.2,1.15)node {$v_6$};

\draw[fill](1,1)circle [radius=0.05];
\draw[fill](2,1)circle [radius=0.05];
\draw[fill](3,1)circle [radius=0.05];
\draw[fill](4,1)circle [radius=0.05];
\draw[fill](5,1)circle [radius=0.05];
\draw[fill](6,1)circle [radius=0.05];

\draw(2.5,0)--(1,1)--(1.2,1.8)--(2,1);
\draw(4.5,0)--(1,1);
\draw[fill](1.2,1.8)circle [radius=0.05];
\draw(0.9,2.1)node {$a_{11}$};
\draw(1,1)--(1.8,1.8)--(2,1);
\draw[fill](1.8,1.8)circle [radius=0.05];
\draw(2.1,2.1)node {$a_{12}$};

\draw(2.5,0)--(6,1);
\draw(2.5,0)--(2,1);
\draw(4.5,0)--(6,1);
\draw(4.5,0)--(2,1);
\draw(3,1)--(3.5,1.8)--(4,1);
\draw[fill](3.5,1.8)circle [radius=0.05];
\draw(3.5,2.1)node {$a_3$};

\draw(2.5,0)--(5,1)--(5.2,1.8)--(6,1);
\draw(4.5,0)--(5,1);
\draw[fill](5.2,1.8)circle [radius=0.05];
\draw(4.8,2.1)node {$a_{21}$};
\draw(5,1)--(5.8,1.8)--(6,1);
\draw[fill](5.8,1.8)circle [radius=0.05];
\draw(6.1,2.1)node {$a_{22}$};

\draw[fill](2.5,0)circle [radius=0.05];
\draw(2.6,-0.4)node {$b_{11}$};
\draw[fill](4.5,0)circle [radius=0.05];
\draw(4.6,-0.4)node {$b_{12}$};

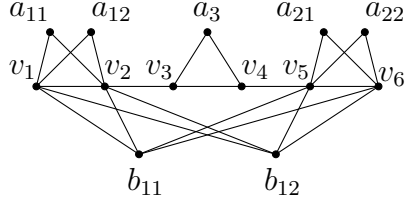
\draw(3.7,-1.5)node{\scriptsize Figure 4.3. The mcc-whiskered graph $G^{mc}$ based on $e_{\pi_2}$};
\end{tikzpicture}
\end{center}

Note that neither  $G^{cc}$ nor $G^{mc}$ is unmixed in the example. Note also that the cardinalities  of  the maximal independent sets of $G^{cc}$  are between 6 and $6+2$, where  $\{v_1,v_2,a_i\mid 3\le i\le 6\}$, $\{b_1, v_2, v_4,v_6,a_1,a_3,a_5\}$ and $\{a_i,b_j\mid 1\le i\le 6,1\le j\le 2\}$ are maximal independent vertex subsets of the graph $G^{cc}$, of cardinality 6,7,8 respectively.

}
\end{myeg}

Recall that a simplicial complex $\D$ is called {\it $(m,n)$ - pure}, if for each facet $F$ of $\D$, $m\le dim\, F\leq n$ holds.

We begin with  the following immediate observation:

\begin{myprn} \label{mccpure} For any graph $G$, let $\pi=\{W_1,\ldots,W_d\}$ be a clique-partition of the vertex set $V(G)$, let
$e_\pi=\{U_1,\ldots, U_s\}$ be a clique cluster decomposition based on $\pi$. Assume further that in the set $e_\pi$, there are totally $s-r$ clusters $U_i$ (of cliques) that is identical with some $W_i$. Let $\D={\rm Ind}\,G^{cc}$ be the independent simplicial complex of the cc-whiskered graph $G^{cc}$. Then the simplicial complex $\D$ is $(d-1,d+r-1)$ - pure of dimension $d+r-1$.
\end{myprn}

\begin{prof} Assume that $W=:\{a_1,\ldots,a_d;b_1,\ldots,b_r\}$ are the added vertex set in $G^{cc}$, then $W$ is a maximal independent vertex set of the graph $G^{cc}$. Note that each subgraph  of the graph $G^{cc}$ induced on $W_i\cup\{a_i\}$ is a clique ($1\le i\le d$), thus for each maximal independent set $V$ of $G^{cc}$, we have $|V\cap (W_i\cup\{a_i\})|\leq 1$, hence $|V|\leq d+r$. This shows ${\rm dim} ({\rm Ind}\,G^{cc})=d+r-1$.

 Assume further $U_i=W_{k_i}$ ($r+1\le i\le s$), then for $r+1\le i\le s$, each vertex of $U_i$ is adjacent  to exactly one vertex of $W$, while each vertex of $U_i$ ($1\le i\le r$) is adjacent  to exactly two vertices of $W$. Now assume that $A=\{a_1,\ldots, a_d\}\seq V(G), B=\{w_1,\ldots,w_r\}\seq V(G)$ be such that $\{\{a_i,v_i\},\{b_j,w_j\}\}\seq E(G^{cc})$. Then for any maximal independent vertex set $V$ of $G^{cc}$, assume further that
 $$\{a_i,\ldots,a_u,b_1,\ldots,b_v \}\seq V, \, V\cap\{a_{u+1},\ldots,a_d,b_{v_1},\ldots,b_r\}=\eset.$$
 Then  it follows from $|S\cup T|=|S|+|T|-|S\cap T|\ge |S|$ and the facts $V\cap \{a_i,v_i\}\not=\eset, V\cap \{b_i,w_i\}\not=\eset $ that $|V|\ge d$. This shows that the simplicial complex ${\rm Ind}\,G^{cc}$ is $(d-1,d+r-1)$ - pure of dimension $d+r-1$.   \end{prof}

\vs{3mm} Note that the cardinalities of maximal independent facets of ${\rm Ind}\,G^{cc}$ may fill the interval $[d,d+r]$, as Example \ref{mccExample} shows.

 Then we proceed with the following observation:

\begin{mylem} \label{mcccases} For a graph $G$, let $\pi$ be a \vcp\,of $G$, $e_\pi$ a clique cluster decomposition based on $\pi$. Let  $G^{md}$ be a mcc-whiskered graph \,of $G$, and let $v$ be any vertex of $G$. Then

$(1)$ There exists a \vcp\,\,$\pi_1$ of $V(G\smallsetminus v)$ and some clique cluster decomposition $e_{\pi_1}$, such that $G^{md}  \smallsetminus v$ is equal to some mcc-whisked graph
$(G\smallsetminus v)^{mc}$ based on $e_{\pi_1}$, together with a possible isolated \vd\, graph;

$(2)$ There exists a \vcp\,\,$\pi_2$ of $V(G\smallsetminus N_G[v])$ and some clique cluster decomposition $e_{\pi_2}$, such that $G^{mc}\smallsetminus N_{G^{md}}[v]$ is equal to some mcc-whiskered graph $(G\smallsetminus N_G[v])^{md}$ based on $e_{\pi_2}$, together with a possible isolated \vd\, graph.

\end{mylem}

\begin{prof} Let $\pi=\{W_1,\ldots,W_d\}$ be a vertex clique-partition of a graph $G$; let $e_\pi=\{U_1,\ldots,U_s\}$
 be a vertex clique cluster-partition based on $\pi$. Assume that $G^{md}$ be some mcc-whiskered graph of type $(d,r)$, constructed as above after adding nonempty finite sets $$A_1,\ldots,A_d;\,\,B_1,\ldots,B_r.$$ For any $1\le i\le d$, take a vertex $v\in W_i\seq V(G)$.

$(1)$ If $W_i=U_{k_i}$ for some $k_i$, we have two subcases. The first case is $|U_{k_i}|>1$; in the case, clearly $G^{md}\smallsetminus v=(G\smallsetminus v)^{md}$ and is a mcc-whiskered graph of type $(d,r)$ on the graph $G\smallsetminus v$;  the other case is $|U_{k_i}|=1$, then the graph $G^{md}\smallsetminus v$ is a mcc-whiskered graph of type $(d-1,r)$ on the graph $G\smallsetminus v$, together with the \vd\, graph $A_{k_i}$, in which $E(G^{md}\smallsetminus v, A_{k_i})=\eset$.

Now assume that $W_i$ is contained in $U_1$ and, $U_1$ is a union of $u$ mutually independent cliques $W_j$'s. If further either $u>2$ or, $u=2$ and $|W_i|\ge 2,$ then $G^{md}\smallsetminus v$ is a mcc-whiskered graph of the type $(d,r)$ on the graph $G\smallsetminus v$; if $u=2$ and $|W_i|=1$, then the graph $G^{md}\smallsetminus v$ is a mcc-whiskered graph of type $(d-1, r-1)$, on the graph $G\smallsetminus v$.

In conclusion, $G^{md}\smallsetminus v$ is a mcc-whiskered graph on the graph $G\smallsetminus v$, together with some possible isolated vertices.

$(2)$ The graph $G^{md}\smallsetminus N_{G^{mc}}[v]$  seems to be   a little bit more complicated, but in fact, only one situation occurs. Assume $v\in W_i\seq U_{k_i}$; note that $U_{k_i}\cup A_i\cup B_{k_i}\seq N_{G^{mc}}[v]$ holds true, and $v$ is not adjacent to other vertices of $W_l\seq U_{k_i}$ if such $l$ with $l\not=k_i$ exists, thus $G^{md}\smallsetminus N_{G^{md}}[v]$ is a mcc-whiskered graph of the type
 $(d-u+t,r-t-w)$ on the base graph $G\smallsetminus N_G[v]$, together   with some possible isolated  \vd\ graphs, where $u$ is the number of $W_j$'s with $U_j\seq N_{G^{md}}[v] $, $t $ is the number
of $U_j$'s which consists of two $W_k$ and exactly one $W_k$
is contained in $N_{G^{md}}[v]$ and, $w$ is the  number  of $U_k$'s with $U_k\seq N_{G^{mc}}[v]$.

This complete the proof.
\end{prof}

\vs{3mm}Note that there is an analogous result if all $md$ are replaced by $mc$, where the phrase “together with a possible isolated \vd\, graph” is replaced by “together with some isolated vertices”.

\vs{3mm}The following known result is also needed:

\begin{mylem}\label{component} For any graph $G$, $G$ is vertex decomposable if and only if all connected components of $G$ are vertex decomposable.
\end{mylem}

For a given graph $G$, the choice of clique partitions $\pi$ are generally varied and abundant; for any vertex clique partition $\pi$ of $G$, the choice of clique cluster partitions $e_\pi$ are   generally also varied and abundant.
Thus starting from any graph $G$, a lot of distinct graphs $G^{mc}$ may be constructed; yet, these graphs share a common good property, i.e., they are \vd, as is shown
in the following main result of this section; note that when passing from $G^{mc}$ to $(G\smallsetminus v)^{mc}$ (and, to $(G\smallsetminus N_G[v])^{mc}$), both the related $\pi$ and $e_\pi$ have changed somewhat.

\begin{mthm} \label{mccPVD} For any graph $G$, let $\pi$ be a \vcp\,of $G$, $e_\pi$ a clique cluster decomposition based on $\pi$ and, $G^{md}$ a  multiple clique cluster-whiskered graph on the base graph $G$. Then the independence simplicial complex ${\rm Ind}\,G^{md}$ of the graph $G^{md}$ is nonpure vertex decomposable.
\end{mthm}

\begin{prof} We use induction on $|V(G)|$. If $G$ has only one end vertex, the result is clear. In the following, assume $|V(G)|>1$ and assume that the conclusion holds true for graphs whose vertex set has less number of vertices. Then for any vertex $v\in W_1$ of $G$, by Lemma \ref{mcccases}, $G^{md}  \smallsetminus v$ is a mcc-whiskered graph together with a possible isolated \vd\,graph $B_1$ added, thus $G^{md}  \smallsetminus v$ is \vd\, by induction hypothesis and Lemma 4.4; $G^{md}\smallsetminus N_{G^{md}}[v]$ is a mcc-whiskered graph together with some possible isolated \vd\,graphs, thus it is also \vd. Furthermore, if $W$ is an independent vertex set of the graph $G^{md}\smallsetminus N_{G^{md}}[v]$, then any $a_1\in A_1$ is in $N_{G^{md}}[v]$ such that $W\cup \{v\}$ is a larger independent vertex set of the graph $G^{md}  \smallsetminus v$. Thus by Definition \ref{VD graphversion}, each vertex of $G$ is a shedding vertex of $G^{md} $; in particular, the graph $G^{md} $ is vertex decomposable, i.e., the simplicial complex ${\rm Ind}\,G^{md}$ is vertex decomposable.
\end{prof}

\begin{mycor} For any graph $G$ and a vertex clique-partition $\pi$, let $e_\pi$ be a clique cluster-partition based on $\pi$ and, $G^{md}$ a mcc-whiskered graph on the graph $G$. Then the edge ideal $I(G^{md})$  is sequentially Cohen-Macaulay.

\end{mycor}

\vs{3mm}Note that each vertex of $G$ is a shedding vertex of $G^{mc}$; note that $G^{mc}$ is not unmixed in general. For a squarefree monomial ideal $I$ of a polynomial ring $S$, recall from \cite{Herzog and Hibi 1999} that $S/I$ is sequentially Cohen-Macaulay if and only if the Alexander dual $I^\vee$ of $I$ is componentwise linear.

By the above mentioned constructions, we immediately observe the following

\begin{myprn}\label{mccchordal}
Let $G$ be a simple graph, and let $G^{cc}$ be some clique cluster-whiskered graph  with respect to a clique cluster partition $e_\pi$  based on a \vcp ~$\pi$ of $G$. Then the graph $G^{cc}$ is chordal if and only if $G$ is chordal.
\end{myprn}

Recall that chordal graphs are an important class of vertex decomposable graphs (\cite{Herzog Hibi and Zheng, FT2007, Villareal,Woodroofe}), while a matroid  is vertex decomposable at each vertex.
If we compare Proposition \ref{mccchordal} with Theorems \ref{mccPVD}, it is clear that Theorem \ref{mccPVD}  provides  ways for obtaining large class of  nonpure vertex decomposable graphs that is not chordal. Note that in \cite{Villareal}, among the complete list of 29 connected Cohen-Macaulay graphs {\it with at most six vertices}, 20 are of the form  $G^\pi$, and 17 are of the form $G^{WE}$ (i.e., obtained by adding some whiskers and (or) ears), for some graphs $G$.


In order to count the graded Betti numbers of the \mcwg\, $G^{mc}$, sometimes we may need the following

\begin{mylem} {\rm(\cite[Corollary 2.2]{Jacques2005})}\label{JAC} Let $\D_1$ and $\D_2$ be simplicial complexes with $V(\D_1)\cap V(\D_2)=\eset$. Then
$$\b_{i,d}(\D_1*\D_2)=\sum_{p+q=i}\sum_{r+s=d}\b_{p,r}(\D_1)\b_{q,s}(\D_2).$$
\end{mylem}

For a simplicial complex $\D$ with $n$ vertices and a field $k$, let $S=k[x_1,\ldots,x_n]$ be a polynomial ring with $n$ variables. Recall that the graded Betti numbers are defined by $\b_{i,j}^k(\D)=\b_{i,j}(S/I_\D)$, and that $\b_{0,0}(\D)=1$ and $\b_{0,j}=0$ holds for all $j>0$, thus we only consider $\b_{i,j}(\D)$ with $i>0$, where $\b_{i,j}^k(\D)$ is abbreviated as $\b_{i,j}(\D)$. Note that by definition, we have $$\b_{j,k}(S/I)=\b_{j-1,k}(I),\, \forall k\ge 0,\,\forall j>0.$$

Recall the following recent result on vertex decomposable complexes:


For a graph $G$, let $\D={\rm Ind}\,G$. Note that $I_{\D^\vee}$ is identical with $I(\ol{\D})$, whose minimal monomial generating set is $\{X_F\mid F\in \mathcal M(G)\}$, where $\mathcal M(G)$ consists of all minimal vertex covers of the graph $G$. $I(\ol{\D})$ is called the {\it cover ideal} of $G$, and is denoted as $I_c(G)=:I_{\D^\vee}=I(\ol{\D})$.
Now we are ready to use \cite[Corollary 2.11(1)]{Moradi2} (or, \cite[Corollary 2.7]{FHT2009}) and the proof of Theorem \ref{mccPVD} to derive the following:

\begin{mycor}\label{BN}
For any graph $G$,  let $\pi$ be a \vcp\,of $G$, $e_\pi$ a clique cluster partition based on $\pi$ and, let $G^{mc}$ be some multiple clique cluster-whiskered graph on $G$. Let $v$ be a vertex of $G$, $G_1=G\smallsetminus v,\,G_2= G\smallsetminus N_{G}[v]$, and obtain the corresponding clique cluster-whiskered graphs $(G_1)^{mc}$ and $G_2^{mc}$ as in the proof of Lemma \ref{mcccases}. Then we have the following recursive formulae for the cover ideal  $I_c(G^{mc})$ of the graph $G^{mc}$\,:
$$\b_{i,j}(I_c(G^{mc}))=\b_{i,j-1}(I_c(G_1^{mc}))+\b_{i,j}(I_c(G_2^{mc}))+\b_{i-1,j-1}(I_c(G_2^{mc})).$$

\end{mycor}

\begin{prof} In the proof of Lemma \ref{mcccases}, if $\{v\}=W_1$, then there exists an integer $m=|A_1|$ and a large enough $r$, such that the facet set of $({\rm Ind}\,G^{mc})\smallsetminus v$ consists of all $[r+1,r+m]\cup F_1$, where $F_1$ is a facet of ${\rm Ind}\,(G\smallsetminus v)^{mc}$; thus we have
$$\ol{({\rm Ind}\,G^{mc})\smallsetminus v}=\ol{{\rm Ind}\,(G\smallsetminus v)^{mc}},$$
hence $I(\ol{({\rm Ind}\,G^{mc})\smallsetminus v})=I_c((G\smallsetminus v)^{mc})$ holds true. Similarly, for $\D={\rm Ind}\,G^{mc}$, we also have $I(\ol{lk_\D(v)})=I_c((G\smallsetminus N_{G}[v])^{mc})$. Then the formulae follows from \cite[Corollary 2.11(1)]{Moradi2}.
\end{prof}

\vs{3mm}Corollary \ref{BN} actually also implies the following formulae concerning the projective dimension and Castelnuovo-Mumford regularity.

\begin{mycor} For any graph $G$ with $n$ vertices, let $S=k[x_1,\ldots,x_n]$ be a polynomial ring, let $G^{mc}$ be some \mclwg\, on $G$, and let $\D={\rm Ind}\,G^{mc}$. Let $v$ be a vertex of $G$, let $\D_1$ be the simplicial complex ${\rm Ind}\,(G\smallsetminus v)^{mc}$, and let $\D_2$ be the simplicial complex ${\rm Ind}\,(G\smallsetminus N_{G}[v])^{mc}$. Then we have

$(1)$ ${\rm pd}(S/I_\D)={\rm max}\{{\rm pd}(S/I_{\D_2}), \,\,{\rm pd} (S/I_{\D_1})+1\}$.

$(2)$  ${\rm reg}(S/I_\D)={\rm max}\{{\rm reg}(S/I_{\D_2}),\,\, {\rm reg} (S/I_{\D_1})+1\}\}$.

\end{mycor}

Note that for a \cwg \, $G^\pi$ (and its generalization $\D_{\chi}$, in simplicial case), there exists results better than Corollary \ref{BN}, as the following result of \cite[Theorem 13]{BT2013} shows:

{\it  Let $f=(f_{-1},f_0,\ldots,f_d)$ be the $f$-vector of a $d$-dimensional simplicial complex $\D$, and let $\chi$ be an s-colouring of $\D$. Then the graded Betti numbers of {\small $I_{\D_{\chi}^\vee}$} in $S$ are given by the following formula
$$\beta_{i,i+n}(I_{\D_{\chi}^\vee})=\sum_{j=1}^{d+1}\begin{pmatrix}j\\i\end{pmatrix}f_{j-1}(\D).$$
In particular, ${\rm proj-dim}(I_{\D_{\chi}^\vee})={\rm reg}(S/I_{\D_{\chi}^\vee})=d+1.$}

The unmixed property of $G^\pi$ is crucial for this result. Note that $G^{mc}$ is generally not unmixed, and the class of mcc-whiskered graphs are much larger than the class of \cwg s $G^\pi$, thus the recursive formulae given
in Corollary \ref{BN} is applied and effective, as the following example shows:

\begin{myeg}   {\rm Assume $n\ge 1$ and consider the cycle graph $G=C_{4n}$ with  $V(G)=\{v_i\mid 1\le i\le 4n\}$, and the vertices are drawn clockwise. Let the
vertex clique partition be $V(G)=\cup_{i=1}^{2n}W_i$ where $W_i=\{v_{2i-1},v_{2i}\}$; let the clique cluster partition be $V(G)=U_1\cup \cdots\cup U_n$, where $U_1=W_1\cup W_3, U_2=W_2\cup W_4,\cdots,U_n=W_{2n-2}\cup W_{2n}$. Consider the cc-whiskered graph $G^{cc}$ of type $(2n,n)$, we have $G^{cc}\smallsetminus v_{4n}=(G\smallsetminus v_{4n})^{cc}$, which is a cc-whiskered graph of  type $(2n,n)$ on the line graph $L_{4n-1}$; while
$G^{cc}\smallsetminus N_{G_c}[v_{4n}]=(G\smallsetminus N_G[v_{4n}])^{cc},$
which is a cc-whiskered graph of type $(4n-1, n-1)$ on the line graph $L_{4n-2}$. We draw the three cc-whiskered graphs in Figures 4.4, for $n=2$.

\begin{center}
\begin{tikzpicture}[scale=0.9]

\draw(1,-1)--(1,1)--(3,1)--(5,1)--(7,1)--(7,-1);
\draw(1,1)--(2,1.5)--(3,1);
\draw(7,1)--(6,1.5)--(5,1);
\draw(1,1)--(2,0.5)--(3,1);
\draw(7,1)--(6,0.5)--(5,1);

\draw(7,-1)--(2,0.5)--(5,-1);
\draw(1,-1)--(6,0.5)--(3,-1);

\draw[fill](1,1)circle [radius=0.05];
\draw[fill](3,1)circle [radius=0.05];
\draw[fill](5,1)circle [radius=0.05];
\draw[fill](7,1)circle [radius=0.05];
\draw[fill](2,1.5)circle [radius=0.05];
\draw[fill](6,1.5)circle [radius=0.05];
\draw[fill](2,0.5)circle [radius=0.05];
\draw[fill](6,0.5)circle [radius=0.05];

\draw(1,-1)--(3,-1)--(5,-1)--(7,-1);
\draw(1,-1)--(2,-1.5)--(3,-1);
\draw(7,-1)--(6,-1.5)--(5,-1);

\draw[fill](1,-1)circle [radius=0.05];
\draw[fill](3,-1)circle [radius=0.05];
\draw[fill](5,-1)circle [radius=0.05];
\draw[fill](7,-1)circle [radius=0.05];
\draw[fill](2,-1.5)circle [radius=0.05];
\draw[fill](6,-1.5)circle [radius=0.05];

\draw(1,1.25)node {$v_1$};
\draw(3,1.25)node {$v_2$};
\draw(5,1.25)node {$v_3$};
\draw(7,1.25)node {$v_4$};
\draw(1,-1.25)node {$v_8$};
\draw(3,-1.25)node {$v_7$};
\draw(5,-1.25)node {$v_6$};
\draw(7,-1.25)node {$v_5$};

\draw(3.5,-2.2)node{\scriptsize   A cc-whiskered  graph $C_{8}^{cc}$ of type $(4,2)$ };
\end{tikzpicture}
\end{center}

\begin{center}
\begin{tikzpicture}[scale=0.9]

\draw(1,1)--(3,1)--(5,1)--(7,1)--(9,1)--(11,1)--(13,1);
\draw(1,1)--(2,1.5)--(3,1);
\draw(7,1)--(6,1.5)--(5,1);
\draw(9,1)--(10,1.5)--(11,1);
\draw(13,1)--(13,1.5);
\draw(13,1)--(9,0)--(5,1);
\draw(9,0)--(7,1);

\draw(1,1)--(5,0)--(3,1);
\draw(9,1)--(5,0)--(11,1);

\draw[fill](1,1)circle [radius=0.05];
\draw[fill](2,1.5)circle [radius=0.05];
\draw[fill](3,1)circle [radius=0.05];
\draw[fill](5,1)circle [radius=0.05];
\draw[fill](6,1.5)circle [radius=0.05];
\draw[fill](7,1)circle [radius=0.05];
\draw[fill](9,1)circle [radius=0.05];
\draw[fill](10,1.5)circle [radius=0.05];
\draw[fill](11,1)circle [radius=0.05];
\draw[fill](13,1)circle [radius=0.05];
\draw[fill](9,0)circle [radius=0.05];
\draw[fill](13,1.5)circle [radius=0.05];
\draw[fill](5,0)circle [radius=0.05];

\draw(1,1.25)node {$v_1$};
\draw(3,1.25)node {$v_2$};
\draw(5,1.25)node {$v_3$};
\draw(7,1.25)node {$v_4$};
\draw(9,1.25)node {$v_5$};
\draw(11,1.25)node {$v_6$};
\draw(13.3,1.25)node {$v_7$};

\draw(8,-0.8)node{\scriptsize  The delete graph $C_{8}^{cc}\smallsetminus v_{8}=L_7^{cc}$ of type $(4,2)$};
\end{tikzpicture}
\end{center}

\begin{center}
\begin{tikzpicture}[scale=0.9]

\draw(3,1)--(5,1)--(7,1)--(9,1)--(11,1);
\draw(3,1.5)--(3,1);
\draw(7,1)--(6,1.5)--(5,1);
\draw(9,1)--(10,1.5)--(11,1);

\draw(5,0)--(3,1);
\draw(9,1)--(5,0)--(11,1);

\draw[fill](3,1.5)circle [radius=0.05];
\draw[fill](3,1)circle [radius=0.05];
\draw[fill](5,1)circle [radius=0.05];
\draw[fill](6,1.5)circle [radius=0.05];
\draw[fill](7,1)circle [radius=0.05];
\draw[fill](9,1)circle [radius=0.05];
\draw[fill](10,1.5)circle [radius=0.05];
\draw[fill](11,1)circle [radius=0.05];


\draw[fill](5,0)circle [radius=0.05];

\draw(2.7,1.25)node {$v_2$};
\draw(5,1.25)node {$v_3$};
\draw(7,1.25)node {$v_4$};
\draw(9,1.25)node {$v_5$};
\draw(11,1.25)node {$v_6$};

\draw(8,-0.8)node{\scriptsize  The link $C_{8}^{cc}\smallsetminus N_{C_{8}^{cc}}[v_{8}]=L_5^{cc}$ of type $(3,1)$};

\draw(7.5,-1.6)node{\scriptsize  Figure 4.4.};
\end{tikzpicture}
\end{center}

Note that mcc whiskered graph will occur if we continue delete a vertex. By the recursive formulae in Lemma \ref{BN} and mathematical induction,  the calculation of the graded Betti numbers $\b_{i,j}$ of the cover ideal $I_c(C_{n}^{cc})$ are finally reduced to the computation of  $\b_{i,j}(J)$, where $J$ is the ideal of $S=k[x_1,x_2,\ldots, x_r]$ generated by $x_1,x_2\cdots x_r$ ($r\ge 2$). Note that $S/J$ has the following minimal free resolution by CoCoA \cite{CoCoA}:
$$0\longrightarrow S(-r)\longrightarrow S(-1)\oplus S(-(r-1))\longrightarrow S,$$
i.e., $\b_{-1,0}(J)=\b_{1,r}(J)=\b_{0,1}(J)=\b_{0,r-1}(J)=1.$
} \end{myeg}

\end{document}